\documentclass[reqno,centertags]{amsart}
\usepackage{amsmath,amsfonts,amsthm,amssymb}
\usepackage[a4paper]{geometry}
\usepackage{graphicx,color}
\newtheorem*{thm}{Theorem}
\newtheorem*{cor}{Corollary}
\numberwithin{equation}{section}
\begin{document}
\author{Jens Wirth}
\address{Jens Wirth, Department of Mathematics, Imperial College, London SW7 2AZ, UK}
\email{j.wirth@imperial.ac.uk}
\title[Anisotropic thermo-elasticity in 2D]{Anisotropic thermo-elasticity in 2D\\ Part II: Applications}
\date{}

\begin{abstract}
This note deals with concrete applications of the general treatment of anisotropic thermo-elasticity 
developed in the first part [M. Reissig, J. Wirth, {\em Anisotropic thermo-elasticity in 2D - A unified treatment}, Asympt. Anal. ?? (????) ??--??]. We give dispersive decay rates for solutions to the type-1 system of thermo-elasticity for 
certain types of anisotropic media.
\end{abstract}
\maketitle

\section{Elastic operators and thermo-elastic systems}
In the first part, \cite{RW07} we treated the type-1 system of thermo-elasticity for an anisotropic medium
in two space dimensions. This system reads as
\begin{subequations}
\begin{align}
&U_{tt}+A(\mathrm D)U+\gamma\nabla\theta=0\\
&\theta_t-\kappa\Delta\theta+\gamma\nabla\cdot U_t=0
\end{align}
\end{subequations}
for the elastic displacement $U(t,x)\in\mathbb R^2$ and the temperature difference to the equilibrium state $\theta(t,x)\in\mathbb R$. Physical properties of the medium are described by the thermal conductivity $\kappa>0$, the thermo-elastic coupling $\gamma>0$ and the elastic operator $A(\mathrm D)$. Here we assume that its symbol has the structure
\begin{equation}\label{eq:Adef}
   A(\xi) = \mathbb D(\xi)^T \mathcal S\mathbb D(\xi),
\end{equation}
where 
\begin{equation}
  \mathbb D = \begin{pmatrix}\xi_1 & 0 \\ 0 &\xi_2 \\\xi_2&\xi_1 \end{pmatrix}
\end{equation}
is a matrix of first order symbols and
\begin{equation}
   \mathcal S = \begin{pmatrix} \tau_1 & \lambda & \sigma_1\\\lambda & \tau_2 & \sigma_2\\\sigma_1&\sigma_2 & \mu\end{pmatrix}
\end{equation}
contains the elasticity modules of the medium. Usually one makes the assumption that $\mathcal S$ is positive definite such that $A(\mathrm D)$ is a positive and self-adjoint operator. Then the first equation is hyperbolic, while the second one is parabolic. 

In \cite{RW07} the thermo-elastic system was treated under the assumptions

(A1-2)\quad $A(\xi)$ is positive (self-adjoint, 2-homogeneous and real-analytic in $\eta=\xi/|\xi|$);

(A3)\quad $\#\mathrm{spec}\;A(\eta)=2$ for all $\eta\in\mathbb S^1$; 

(A4)\quad $\gamma^2\ne \varkappa_j(\eta)-\mathrm{tr}\;A(\eta)$ for all hyperbolic
directions  (with respect to $\varkappa_j(\eta)\in\mathrm{spec}\;A(\eta)$),

where a direction $\eta\in\mathbb S^1$ is called hyperbolic with respect to $\varkappa_j(\eta)$ if the corresponding (normalised)
 eigenvector $r_j(\eta)$ is perpendicular to the direction $\eta$.  Under these
conditions the system was reformulated as system of first order and Fourier integral representations
in terms of Fourier integral operators with complex phases established. They in turn lead to $L^p$--$L^q$ decay estimates of solutions. 

Here we will focus on these results. It turned out that the decay properties of solutions depend only 
on the neighbourhood of hyperbolic directions. Essential are
the {\em vanishing orders} $\ell_j(\eta)$ of the so-called coupling functions $a_j(\eta)=\eta\cdot r_j(\eta)$ 
for the corresponding eigenvectors and the {\em orders of tangency} $\bar\gamma_j(\eta)$ for the Fresnel curves 
\begin{equation}
S = \{\;  \xi \;|\;1\in\mathrm{spec}\; A(\xi)\;\}=\{\;\varkappa_j(\eta)^{-1/2}\eta\;|\;\eta\in\mathbb S^1,\;\varkappa_j(\eta)\in\mathrm{spec}\;A(\eta)\;\}
\end{equation}
in the hyperbolic direction and for the corresponding sheet. The following theorem collects dispersive decay rates, the statements in \cite{RW07} are more localised in the sense that we decomposed solutions into modes with different decay behaviour micro-locally.

\begin{thm}[Reissig-Wirth, \cite{RW07}]\label{thm}
Let $\bar\eta\in\mathbb S^1$ be a fixed direction, $\chi\in C_0^\infty(\mathbb S^1)$ a cutoff
function for a small neighbourhood of $\bar\eta$ and $\chi(\xi)=\chi(\xi/|\xi|)$. Then the following
micro-local dispersive estimates hold true:
\begin{enumerate}
\item If $\bar\eta$ is hyperbolic with respect to $\varkappa_j(\bar\eta)$ and
$2\ell_j(\bar\eta)\le \bar\gamma_j(\bar\eta)$, then for any $r>2$ we find a constant $C$ such that
\begin{equation}
  \| \chi(\mathrm D) (U_t, \sqrt{A(\mathrm D)}U, \theta)(t,\cdot)\|_\infty
  \le C (1+t)^{-\frac1{2\ell_j(\bar\eta)}} \big(   \| \chi(\mathrm D)( U_t, \sqrt{A(\mathrm D)}U ,\theta)(0,\cdot)\|_{H^{1,r}}
\big).
\end{equation}
\item If $\bar\eta$ is hyperbolic with respect to $\varkappa_j(\bar\eta)$ and
$2\ell_j(\eta)> \bar\gamma_j(\bar\eta)$, then  for any $r>2$ we find a constant $C$ such that
\begin{equation}
  \| \chi(\mathrm D) (U_t, \sqrt{A(\mathrm D)}U, \theta)(t,\cdot)\|_\infty
  \le C (1+t)^{-\frac1{\bar\gamma_j(\bar\eta)}} \big(   \| \chi(\mathrm D)( U_t, \sqrt{A(\mathrm D)}U ,\theta)(0,\cdot)\|_{H^{1,r}}
\big).
\end{equation}
\item If $\bar\eta$ is not hyperbolic, then  for any $r>2$ we find a constant $C$ such that
\begin{equation}
  \| \chi(\mathrm D) (U_t, \sqrt{A(\mathrm D)}U, \theta)(t,\cdot)\|_\infty
  \le C (1+t)^{-1} \big(   \| \chi(\mathrm D)( U_t, \sqrt{A(\mathrm D)}U ,\theta)(0,\cdot)\|_{H^{1,r}}
\big).
\end{equation}
\end{enumerate}
\end{thm}

We want to apply this statement in several concrete situations. We prepare this in Section~\ref{sec2}
by some general remarks related to the choice of \eqref{eq:Adef} and the resulting achievable
decay rates. Later in Sections~\ref{sec3} to~\ref{sec6} we will discuss isotropic, cubic, rhombic and a case of fully anisotropic media. 

\section{Some general statements}\label{sec2}
The Fresnel curve $S$ has the form $S=S_1\cup S_2$ with $S_j = \{ \varkappa_j(\eta)^{-1/2} \eta \;|\;
\eta\in\mathbb S^1\}$ and
due to (A3) we are allowed to label the (real-analytic) eigenvalues in ascending order,  
$\varkappa_1(\eta)<\varkappa_2(\eta)$. 

Because the entries of the matrix $A(\xi)$ are polynomial of degree two we immediately see that the 
associated Fresnel curve $S$ is algebraic of degree four. Hence, any straight line intersects $S$ in at most four points. As pointed out in \cite{Duff} this has important consequences. Any line intersecting the inner sheet $S_2$ has to intersect the outer one $S_1$ in two points. Thus it can intersect $S_2$ at most twice and therefore we see that $S_2$ has to be strictly convex. Thus the order of tangency $\bar\gamma_2(\eta)$ for tangents on $S_2$ at $\varkappa_2(\eta)^{-1/2}\eta$ is two  for any direction $\eta$. Similarly we see that for the outer one $2\le \bar\gamma_1(\eta)\le 4$. 

The hyperbolic directions are determined by the fact that $\eta$ and $\eta^\perp$ are eigenvectors of
the matrix $A(\eta)$, thus they satisfy $\eta^\perp\cdot A(\eta)\eta=0$. This is a polynomial equation of 
order eight. Introducing polar co-ordinates $\eta=(\cos\phi,\sin\phi)^T$ allows to rewrite it as trigonometric polynomial and hyperbolic directions are zeros of 
\begin{equation}\label{eq:HypDirPol}
  4(\sigma_1+\sigma_2)\cos2\phi+4(\sigma_1-\sigma_2)\cos4\phi - 2(\tau_1-\tau_2)\sin2\phi -(\tau_1+\tau_2-2\lambda-4\mu)\sin4\phi.
\end{equation}
It turns out that this polynomial has either four, six or eight roots according to the choice of $\mathcal S$ or it vanishes identically (for isotropic media).

The vanishing orders of the coupling functions are bounded by the vanishing order of \eqref{eq:HypDirPol}: If $a_j(\eta)$ vanishes to order $k$, then $r_j(\eta)=\eta^\perp$ up to order $k$ and
therefore $\eta\cdot A(\eta)\eta^\perp$ vanishes at least to order $k$. But this means that the 
order of zeros to \eqref{eq:HypDirPol} always dominates the vanishing order of the coupling functions and therefore (provided that the medium is anisotropic)
\begin{equation}\label{eq:ellbound}
   4\le \sum_{j=1,2} \sum_{\eta: a_j(\eta)=0} \ell_j(\eta) \le 8
\end{equation}
for the vanishing orders $\ell_j(\eta)$ of $a_j$ in $\eta$.

These general observations imply that Theorem~\ref{thm} gives as overall decay rate for the medium
either $t^{-1/2}$, $t^{-1/3}$ or $t^{-1/4}$. We will provide examples for all these decay rates.

\section{Isotropic media}\label{sec3}
For isotropic media we have $\sigma_1=\sigma_2=0$ and $\tau_1=\tau_2=\lambda+2\mu$, i.e.
\begin{equation}
  A(\eta) = \begin{pmatrix} (\lambda+\mu)\eta_1^2+\mu & (\lambda+\mu)\eta_1\eta_2 \\
  (\lambda+\mu)\eta_1\eta_2 & (\lambda+\mu)\eta_2^2+\mu \end{pmatrix} = \mu \mathrm{I} + (\lambda+\mu)\eta\otimes\eta.
\end{equation}
In this case
$\mathrm{spec}\;A(\eta)=\{\lambda+2\mu,\,\mu\}$ with corresponding eigenvectors $\eta$ and $\eta^\perp$. Thus one of the coupling functions (the one corresponding to $\mu$) vanishes identically and
the Fresnel curves are two concentric circles with tangency indices $\bar\gamma_1(\eta)=\bar\gamma_2(\eta)=2$ for all $\eta$. Our assumptions read as follows:

(A1-2)\quad $\mu>0$, $\lambda>-2\mu$;

(A3)\quad $\lambda+\mu\ne0$;

(A4)\quad $\gamma^2\ne-(\lambda+\mu)$.

Theorem 3.1 from \cite{RW07} implies the well-known fact that the solution to the thermo-elastic system can be decomposed into three components, one related to the parabolic eigenvalue $\lambda+2\mu$
and two related to the hyperbolic one. The hyperbolic components satisfy wave equations and the
corresponding dispersive estimate with $L^1$--$L^\infty$ decay rate $t^{-1/2}$. The parabolic component decays like $t^{-1}$. The decomposition is just the usual Helmholtz decomposition.

\begin{cor} Isotropic media have the dispersive decay rate $t^{-1/2}$. \end{cor}

\section{Cubic media}\label{sec4}
Cubic media are defined by $\tau_1=\tau_2$ and $\sigma_1=\sigma_2=0$, i.e.
\begin{equation}
  A(\eta) = \begin{pmatrix} (\tau-\mu)\eta_1^2+\mu & (\lambda+\mu)\eta_1\eta_2 \\
  (\lambda+\mu)\eta_1\eta_2 & (\tau-\mu)\eta_2^2+\mu \end{pmatrix}.
\end{equation}
This matrix is positive if and only if 
\begin{equation}
  \mathrm{tr}\;A(\eta)=\tau+\mu > 0
\end{equation}
and
\begin{equation}
  \det A(\eta) = \big((\tau-\mu)^2-(\lambda+\mu)^2\big)\eta_1^2\eta_2^2 +\mu(\tau-\mu)+\mu^2 >0
\end{equation}
for all $\eta\in\mathbb S^1$. Using $0\le \eta_1^2\eta_2^2=\cos^2\phi\sin^2\phi\le\frac14$ we see
that this is equivalent to $\mu\tau>0$ and $4\mu\tau+(\tau-\mu)^2-(\lambda+\mu)^2 =\tau^2+2\mu\tau-\lambda^2-2\mu\lambda=(\tau-\lambda)(\tau+\lambda+2\mu)>0$.

Condition (A3) means that the eigenvalues of $A(\eta)$ are distinct. On the contrary, $A(\eta)$ has multiple eigenvalues for a direction $\eta$ if and only if $A(\eta)$ is a multiple of the identity or
\begin{equation}
   (\lambda+\mu)\eta_1\eta_2=0,\qquad\text{and}\qquad (\tau-\mu)\eta_1^2=(\tau-\mu)\eta_2^2.
\end{equation}
The first equation means $(\lambda+\mu)=0$ or $\eta_1\eta_2=0$. In the latter case the second equation is satisfied for $(\tau-\mu)=0$, in the first one for $\eta_1^2=\eta_2^2$. These conditions describe degenerate directions and (A3) holds if and only if $(\lambda+\mu)(\tau-\mu)\ne0$.

In order to find all hyperbolic directions we have to solve
\begin{equation}
 (\tau-\lambda-2\mu)\sin4\phi=0.
\end{equation}
The case $\tau=\lambda+2\mu$ gives isotropic media and we can skip it. In the remaining cases we require $\sin4\phi=0$, i.e. $\phi=k\pi/4$. This corresponds to eight hyperbolic directions and using the general statements we immediately know that the coupling functions vanish to first order in these directions.

It remains to check condition (A4). We distinguish between the directions $\phi=0$, i.e. $\eta=(1,0)^T$, and $\phi=\pi/4$, i.e. $\eta=(\sqrt2/2,\sqrt2/2)^T$. In the first case $A(\eta)=\mathrm{diag}(\tau,\mu)$ and the hyperbolic eigenvalue (i.e. to the eigenvector $\eta^\perp=(0,1)$) is $\mu$. Thus we require
$\gamma^2\ne\mu-\tau$. For the second hyperbolic direction we obtain similarly the eigenvalues
$\mathrm{spec}\;A(\eta)=\{\mu+(\lambda+\tau)/2,\, (\tau-\lambda)/2\}$. Again the second one corresponds to $\eta^\perp$ and is therefore the hyperbolic eigenvalue. Thus we have to require
$\gamma^2\ne(\tau-\lambda)/2-\mu-(\lambda+\tau)/2=-(\mu+\lambda)$.

We collect our assumptions:

(A1-2)\quad $\mu,\tau>0$, $-\tau-2\mu<\lambda<\tau$;

(A3)\quad $(\lambda+\mu)(\tau-\mu)\ne0$;

(A4)\quad $\gamma^2\not\in\{-(\lambda+\mu),\,\mu-\tau\}$.

\begin{figure}
\includegraphics[width=5cm]{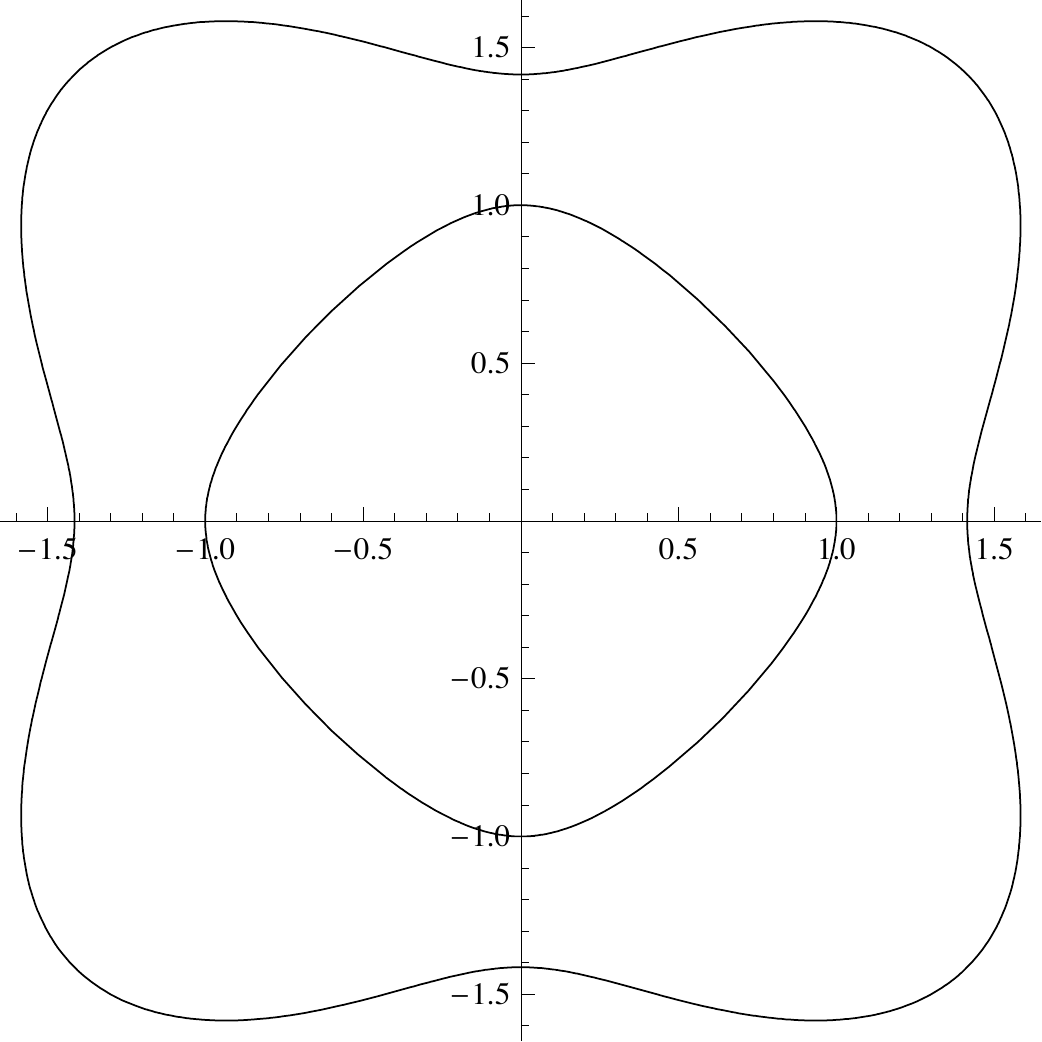} \includegraphics[width=5cm]{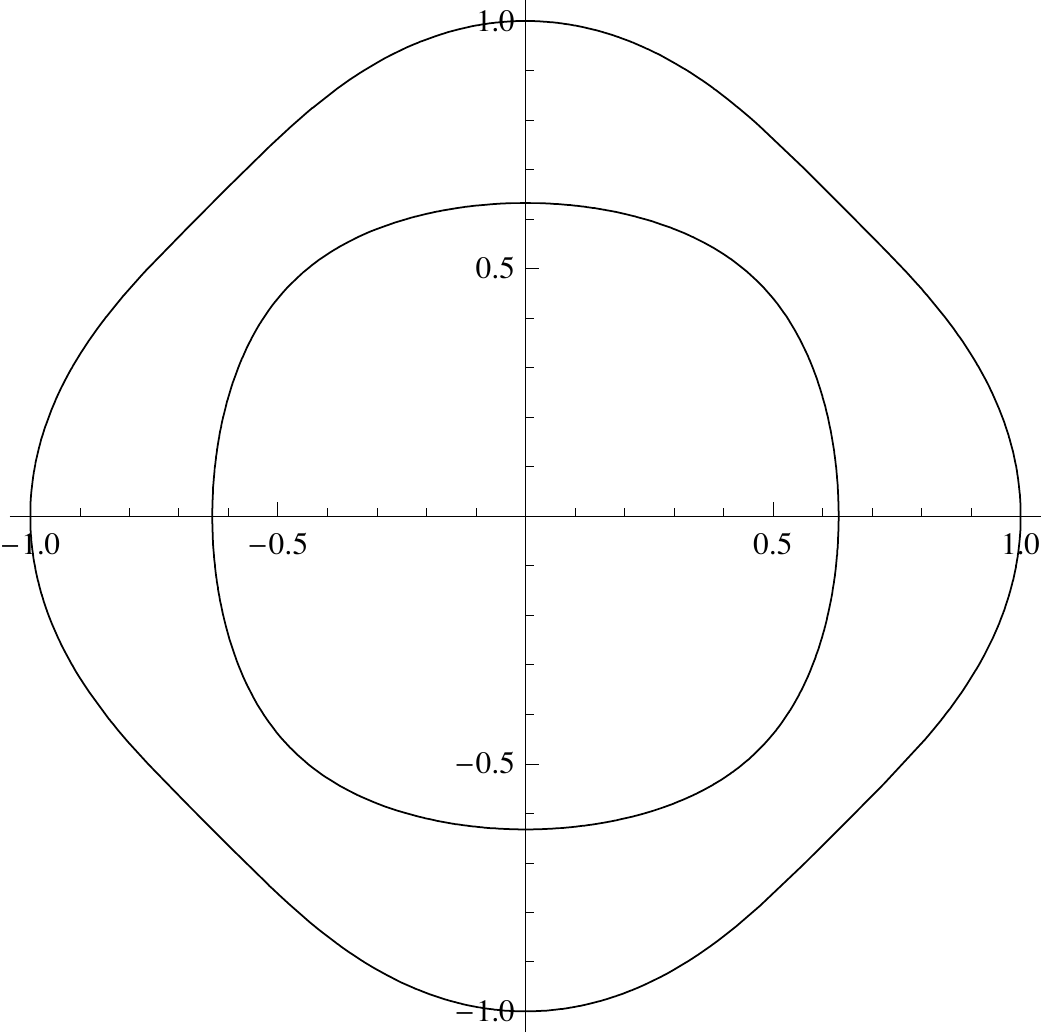}\\
\includegraphics[width=5cm]{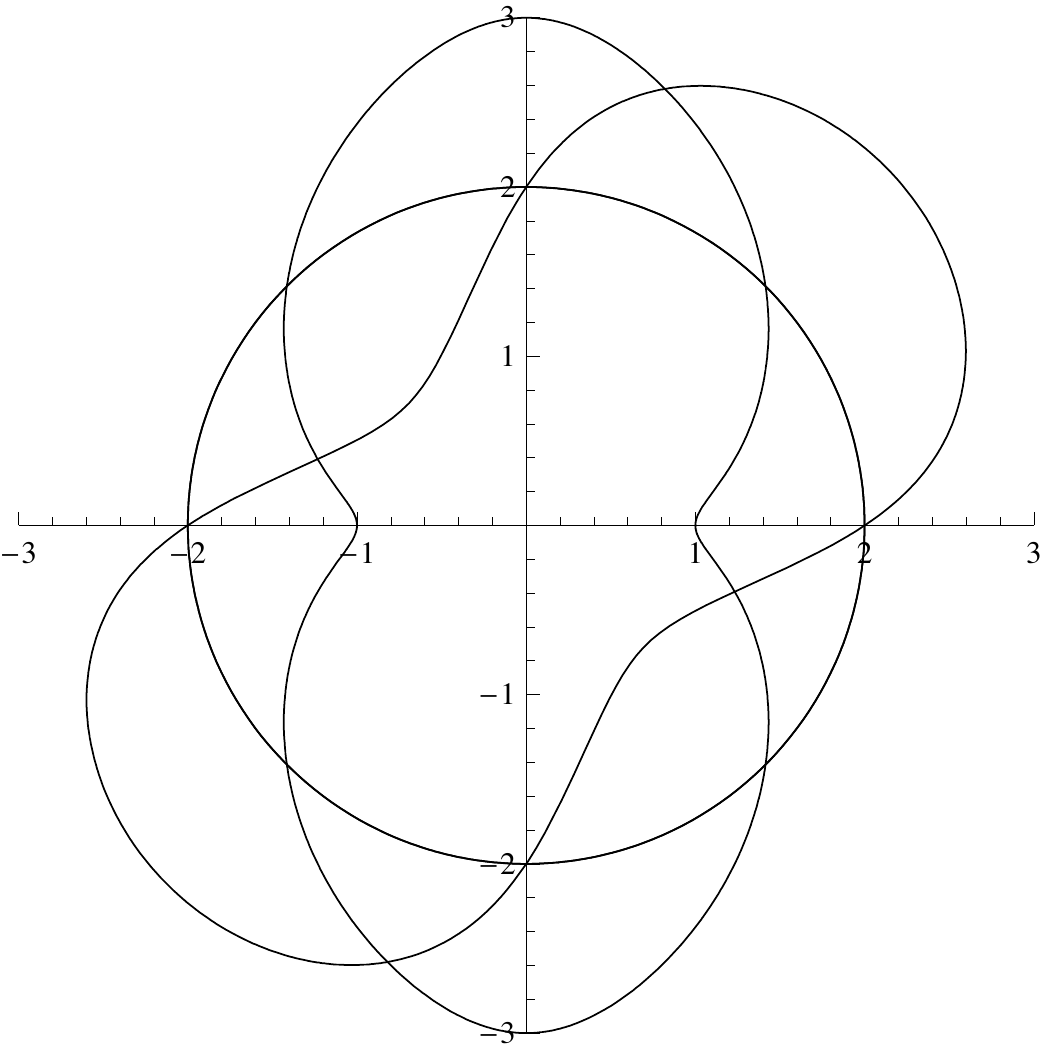} \includegraphics[width=5cm]{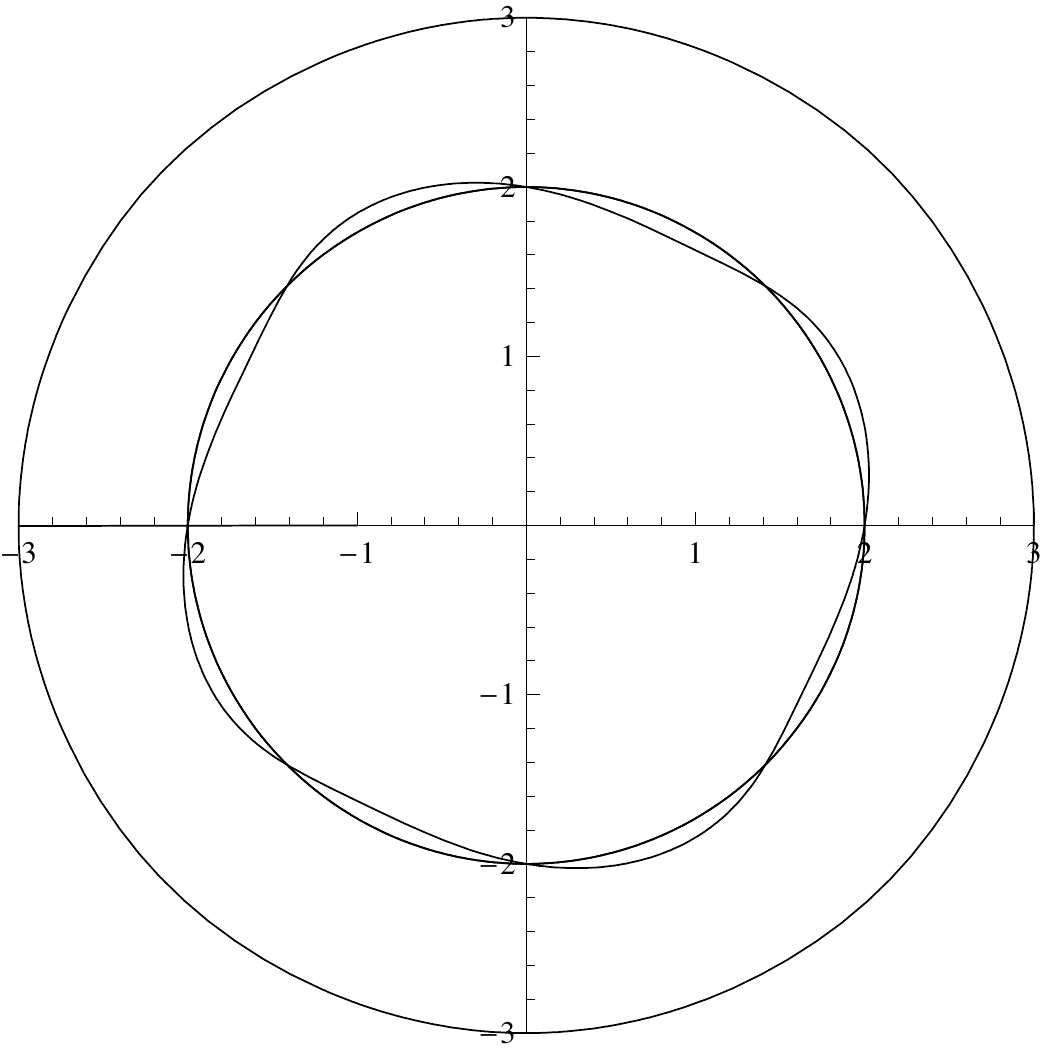}\\[3ex]
\includegraphics[width=5cm]{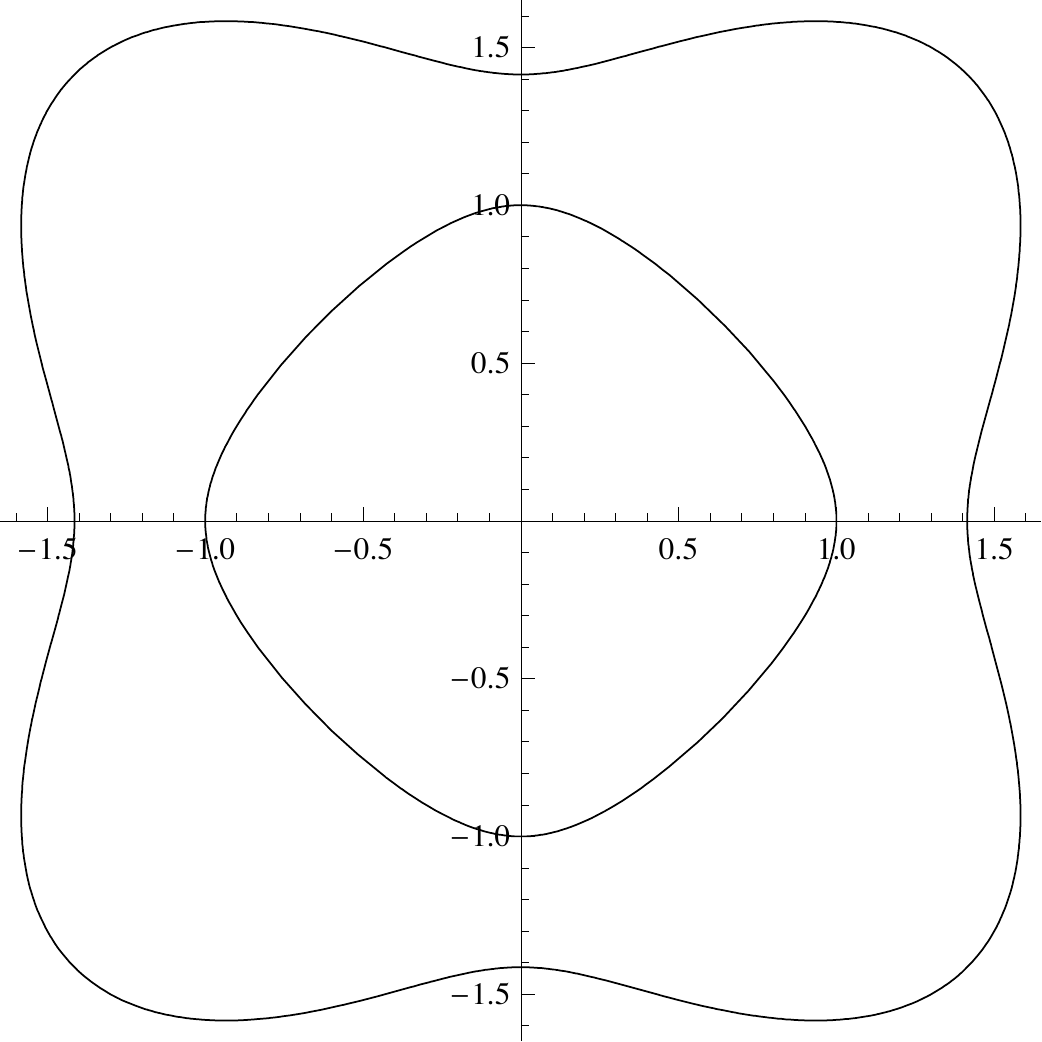} \includegraphics[width=5cm]{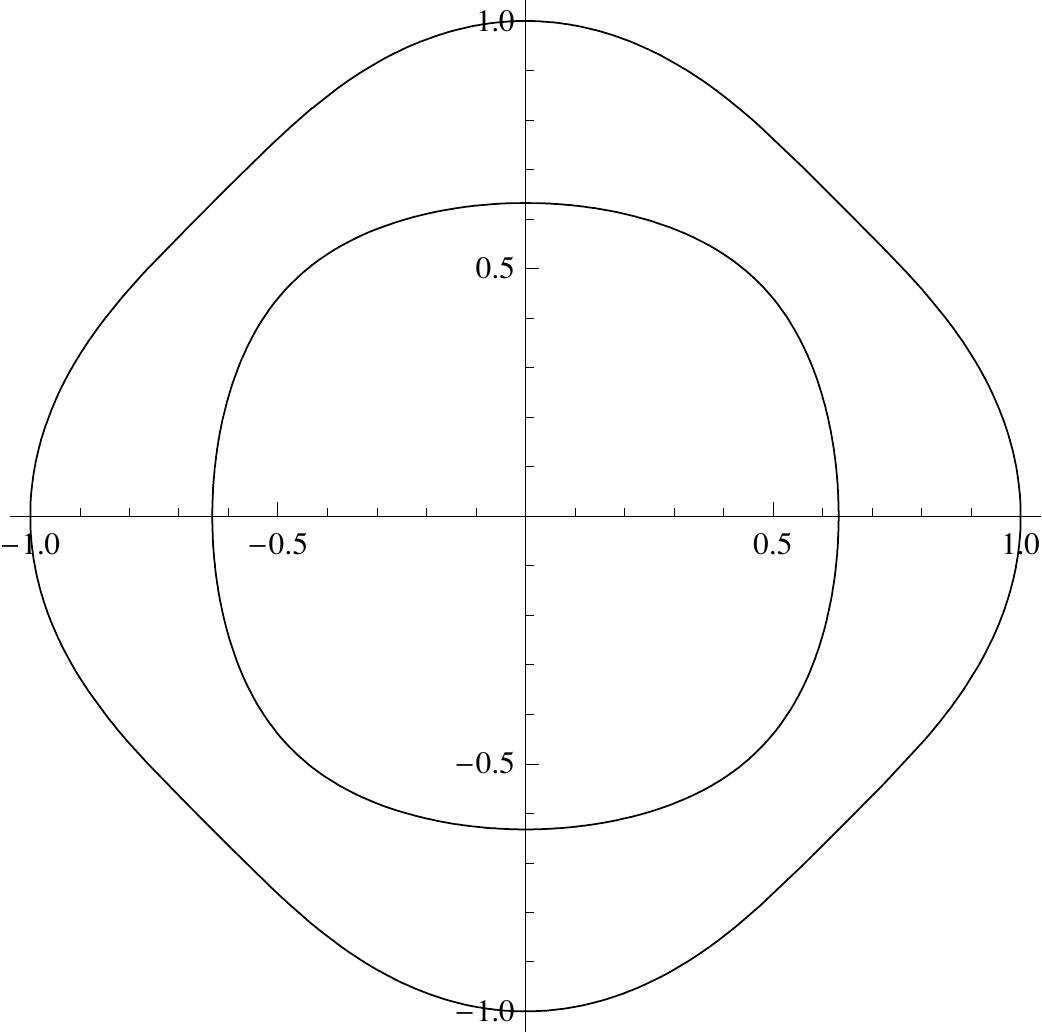}\\
\includegraphics[width=5cm]{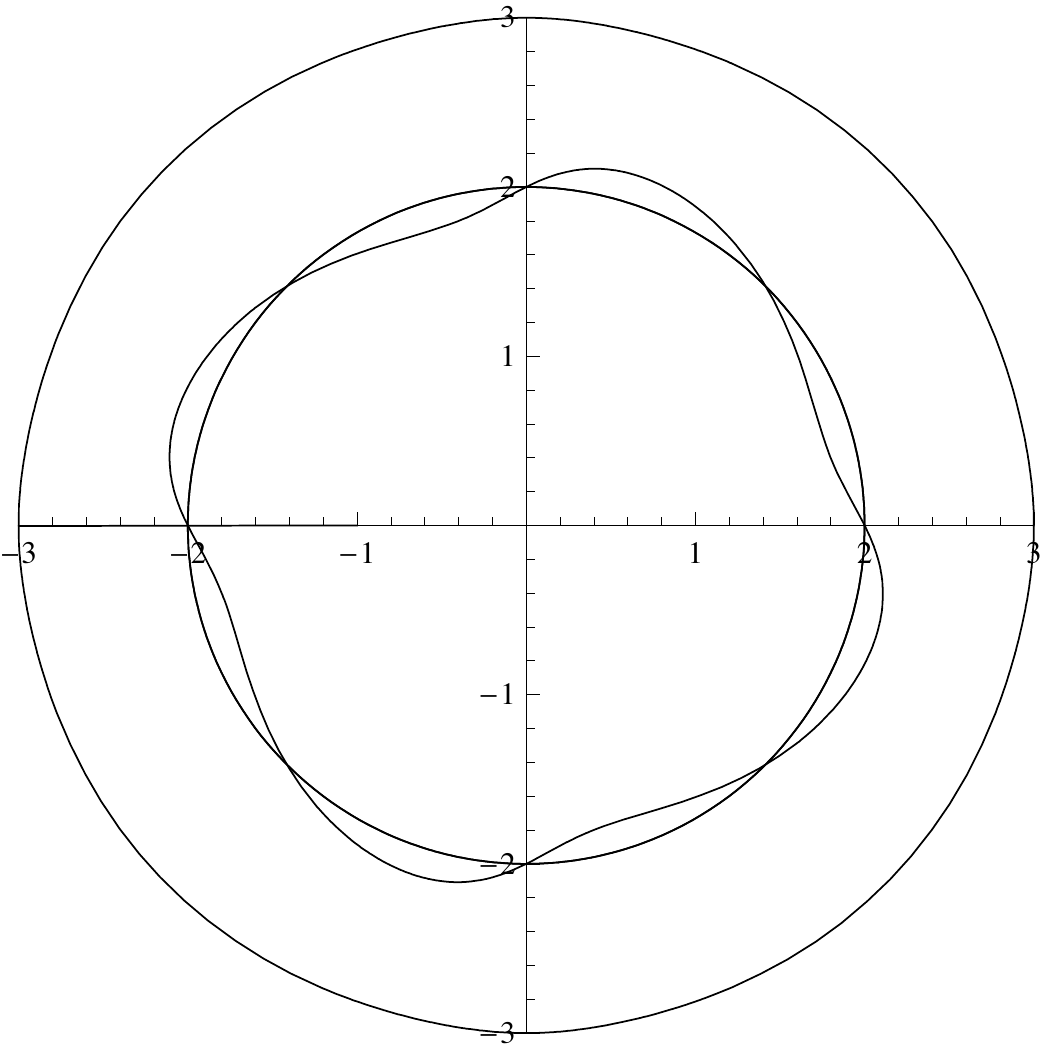} \includegraphics[width=5cm]{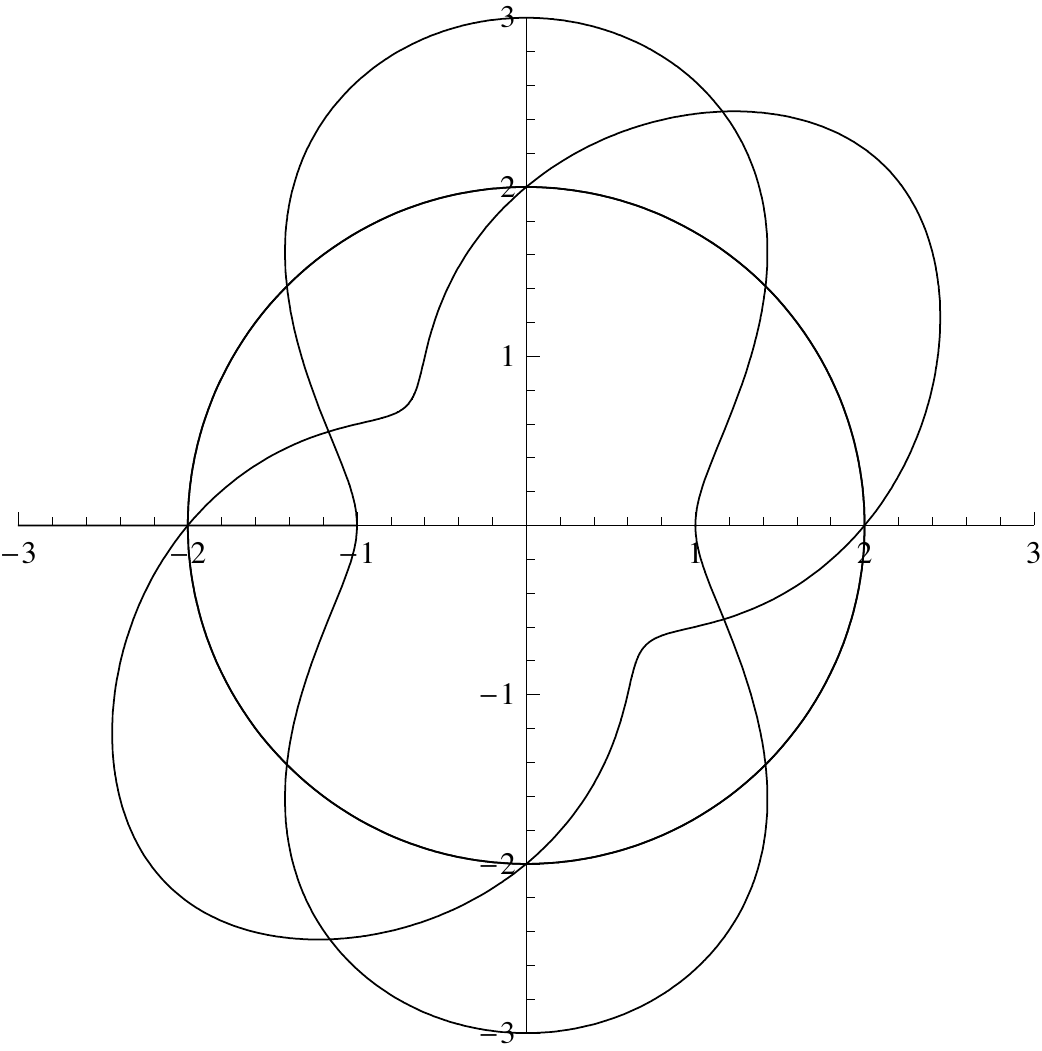}
\caption{Fresnel curves for cubic media together with the corresponding coupling functions (as curves
$\{(2+a_j(\eta))\eta\;|\;\eta\in\mathbb S^1\}$). Depicted are $\lambda=0$, $\mu=1$, $\tau=\frac12$ (upper left), $\lambda=0$, $\mu=1$, $\tau=\frac52$ (upper right),
$\lambda=-2$, $\mu=1$, $\tau=\frac12$ (lower left) and $\lambda=-2$, $\mu=1$, $\tau=\frac52$ (lower right).}
\label{fig1}
\end{figure}

Figure~\ref{fig1} depicts Fresnel curves and coupling functions for different cubic media.

\begin{cor} Cubic media have the dispersive decay rate $t^{-1/2}$. \end{cor}

\section{Rhombic media}\label{sec5}
Rhombic media are determined by the condition $\sigma_1=\sigma_2=0$, i.e. 
\begin{equation}
  A(\eta) = \begin{pmatrix} (\tau_1-\mu)\eta_1^2+\mu & (\lambda+\mu)\eta_1\eta_2 \\
  (\lambda+\mu)\eta_1\eta_2 & (\tau_2-\mu)\eta_2^2+\mu \end{pmatrix}.
\end{equation}
This matrix is positive if  and only if (again using polar co-ordinates)
\begin{equation}
   \mathrm{tr}\;A(\eta) = \mu+\tau_2+(\tau_1-\tau_2)\cos^2\phi>0,
\end{equation}
i.e. if $\mu+\min(\tau_1,\tau_2)>0$, and 
\begin{equation}
  \det A(\eta) = \mu\tau_1\cos^4\phi+\mu\tau_2\sin^4\phi-(\lambda^2+2\lambda\mu-\tau_1\tau_2)\cos^2\phi\sin^2\phi>0,
\end{equation}
i.e. if $\mu\tau_1,\mu\tau_2>0$ and $2\mu\sqrt{\tau_1\tau_2}>\lambda^2+2\lambda\mu-\tau_1\tau_2$.
The last inequality can be simplified to $(\sqrt{\tau_1\tau_2}-\lambda)(\sqrt{\tau_1\tau_2}+\lambda+2\mu)>0$.

Assumption (A3) means that we find no solutions to
\begin{equation}
   (\lambda+\mu)\sin2\phi=0\qquad\text{and}\qquad (\tau_1-\mu)\cos^2\phi=(\tau_2-\mu)\sin^2\phi.
\end{equation}
Thus the first condition gives that (A3) is violated if $\lambda+\mu=0$ or $\phi=k\pi/2$. In the latter case the second
equation gives $\tau_1=\mu$ or $\tau_2=\mu$ (depending on which $\phi$ we chose). If on the contrary
$\lambda+\mu=0$, the second equation reads as $(\tau_1+\tau_2-2\mu)\cos^2\phi=\tau_2-\mu$,
which has further solutions if either $\tau_1,\tau_2>\mu$ or $\tau_1,\tau_2<\mu$. Thus for (A3) we require that $(\tau_1-\mu)(\tau_2-\mu)\ne0$  together with $(\lambda+\mu)\ne0$ if $(\tau_1-\mu)(\tau_2-\mu)>0$. 

Hyperbolic directions are determined by
\begin{equation}
   2(\tau_1-\tau_2)\sin2\phi +(\tau_1+\tau_2-2\lambda-4\mu)\sin4\phi =0.
\end{equation}
This equation has the four zeros $\phi=k\pi/2$ and, dividing by $\sin2\phi$ we see that all further
zeros satisfy $(\tau_1-\tau_2) + (\tau_1+\tau_2-2\lambda-4\mu)\cos2\phi =0$. This condition gives further zeros if $\tau_1-\tau_2\le\tau_1+\tau_2-2\lambda-4\mu$ or $\tau_2-\tau_1\le\tau_1+\tau_2-2\lambda-4\mu$, i.e. if $(\lambda+2\mu-\tau_1)(\lambda+2\mu-\tau_2)\ge0$. If this expression vanishes, but $\tau_1\ne\tau_2$, we get triple zeros for the corresponding $\phi$. If $\tau_1=\lambda+2\mu$, $\tau_2\ne\tau_1$ this holds for $\phi=k\pi$, otherwise for $\phi=k\pi+\frac\pi2$.

\begin{figure}[htbp]
\parbox{5cm}{\includegraphics[width=5cm]{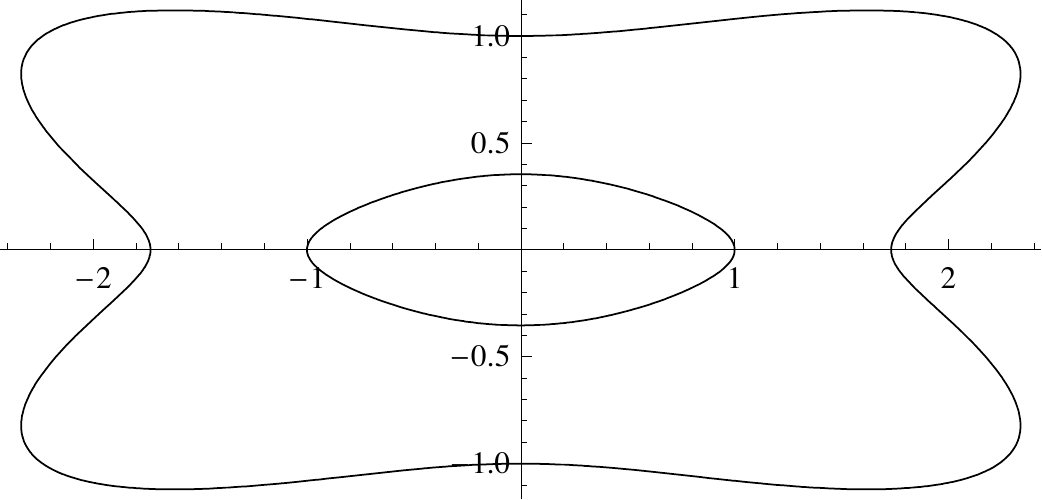}}
\parbox{5cm}{\includegraphics[width=5cm]{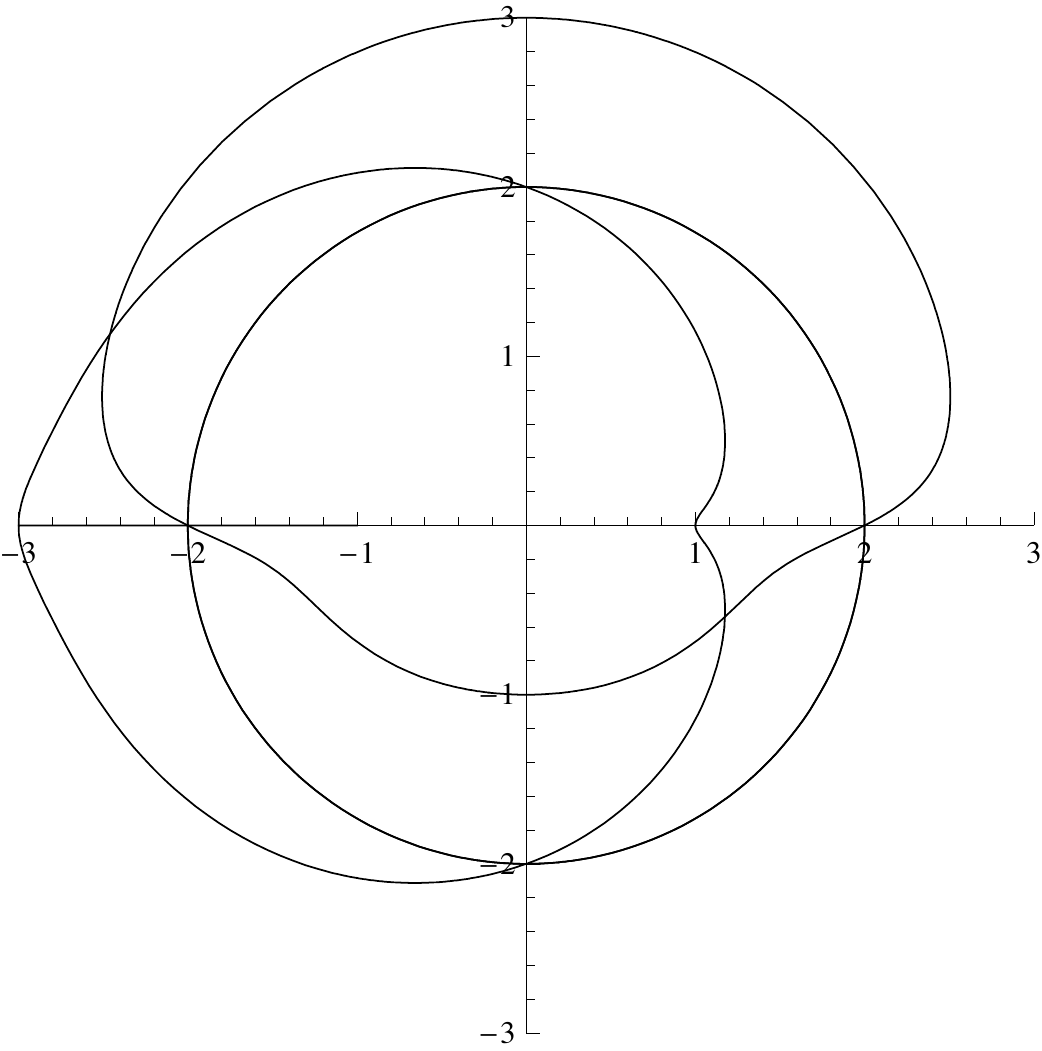}}
\caption{Fresnel curves and coupling functions for a rhombic medium with $\lambda=\frac12$, $\mu=1$,
$\tau_1=4$ and $\tau_2=8$. The medium has eight hyperbolic directions.}
\label{fig2a}
\end{figure}

\begin{figure}[htbp]
\parbox{5cm}{\includegraphics[width=5cm]{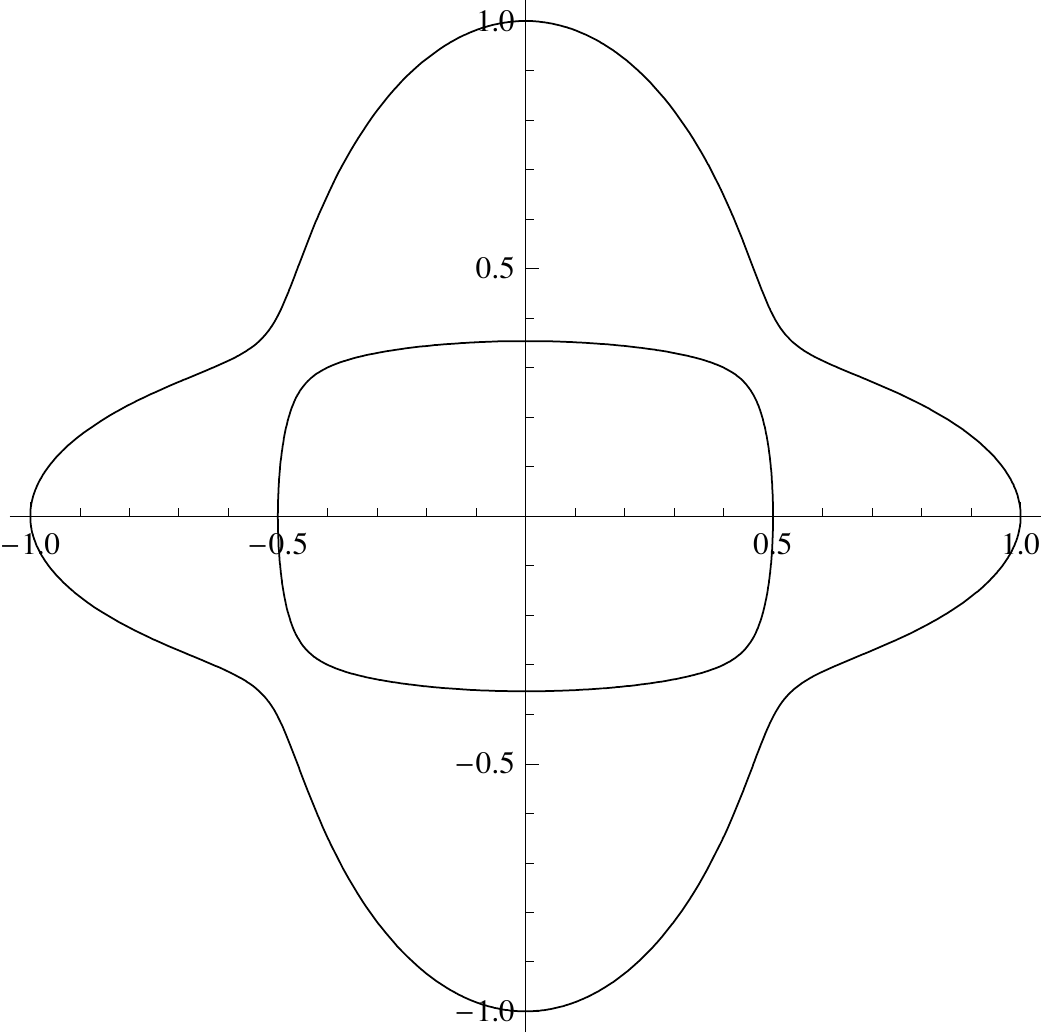}}
\parbox{5cm}{\includegraphics[width=5cm]{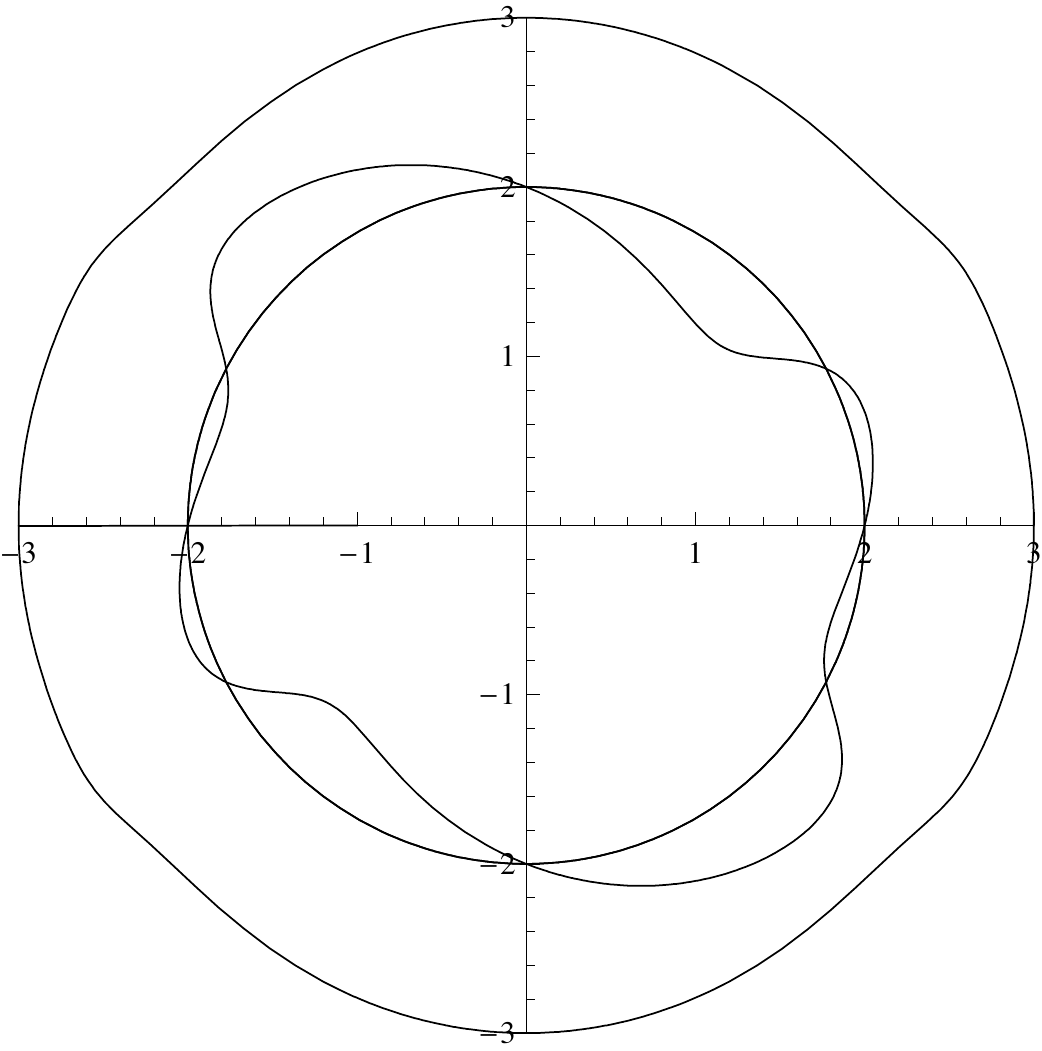}}
\caption{Fresnel curves and coupling functions for a rhombic medium with $\lambda=1$, $\mu=1$,
$\tau_1=\frac13$ and $\tau_2=8$. The medium has four hyperbolic directions.}
\label{fig2b}
\end{figure}

To check condition (A4) we have to determine the corresponding hyperbolic eigenvalues. For 
$\phi=0$ this gives $\gamma^2\ne\mu-\tau_1$, for $\phi=\frac\pi2$ similarly $\gamma^2\ne\mu-\tau_2$.
If the additional direction $\phi=\frac12\arccos (\frac{\tau_1-\tau_2}{\tau_1+\tau_2-2\lambda-4\mu})$ 
exists the difference of the eigenvalues is $\lambda+\mu$ (as for cubic media!) and we require further $\gamma^2\ne-(\lambda+\mu)$. 

In detail: We know that the eigenvectors
of $A(\eta)$ are given by $\eta=(\cos\phi,\sin\phi)$ and $\eta^\perp=(-\sin\phi,\cos\phi)$. This allows 
to determine the eigenvalues by just calculating the products $A(\eta)\eta$ and $A(\eta)\eta^\perp$.
We give only the first components, in the first case it is $(\tau_1-\mu)\cos^3\phi+\mu\cos\phi+(\lambda+\mu)\sin^2\phi\cos\phi$, which gives the (parabolic) eigenvalue $(\tau_1-\mu)\cos^2\phi+(\mu+\lambda)\sin^2\phi+\mu = (\tau_1-2\mu-\lambda)\cos^2\phi+\lambda+2\mu$. For the second one we obtain
$(-(\tau_1-\mu)\cos^2\phi\sin\phi-\mu\sin\phi+(\lambda+\mu)\cos^2\phi\sin\phi$, which by division
through $-\sin\phi$ gives the hyperbolic eigenvalue $(\tau_1-2\mu-\lambda)\cos^2+\mu$. Thus
the difference between both is $-(\lambda+\mu)$.

\begin{figure}[htbp]
\parbox{5cm}{\includegraphics[width=5cm]{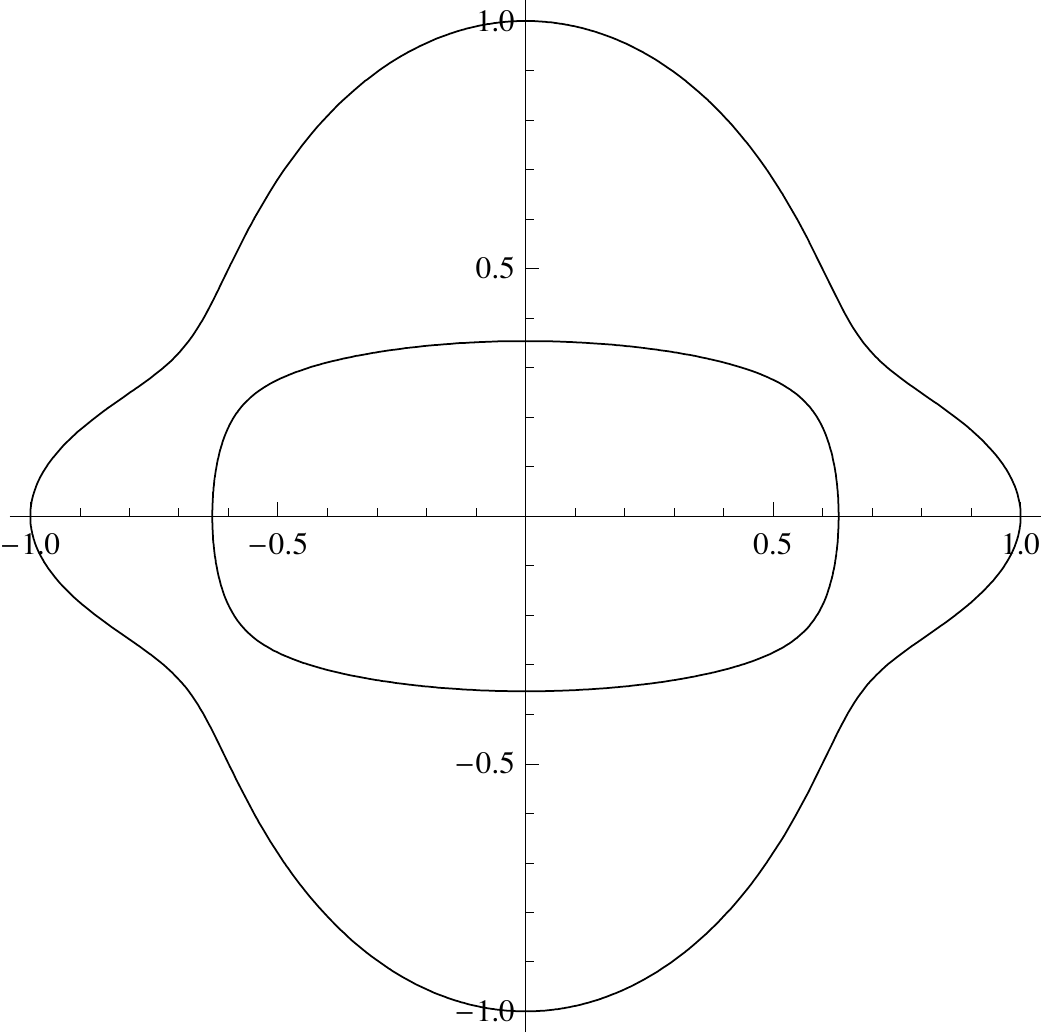}}
\parbox{5cm}{\includegraphics[width=5cm]{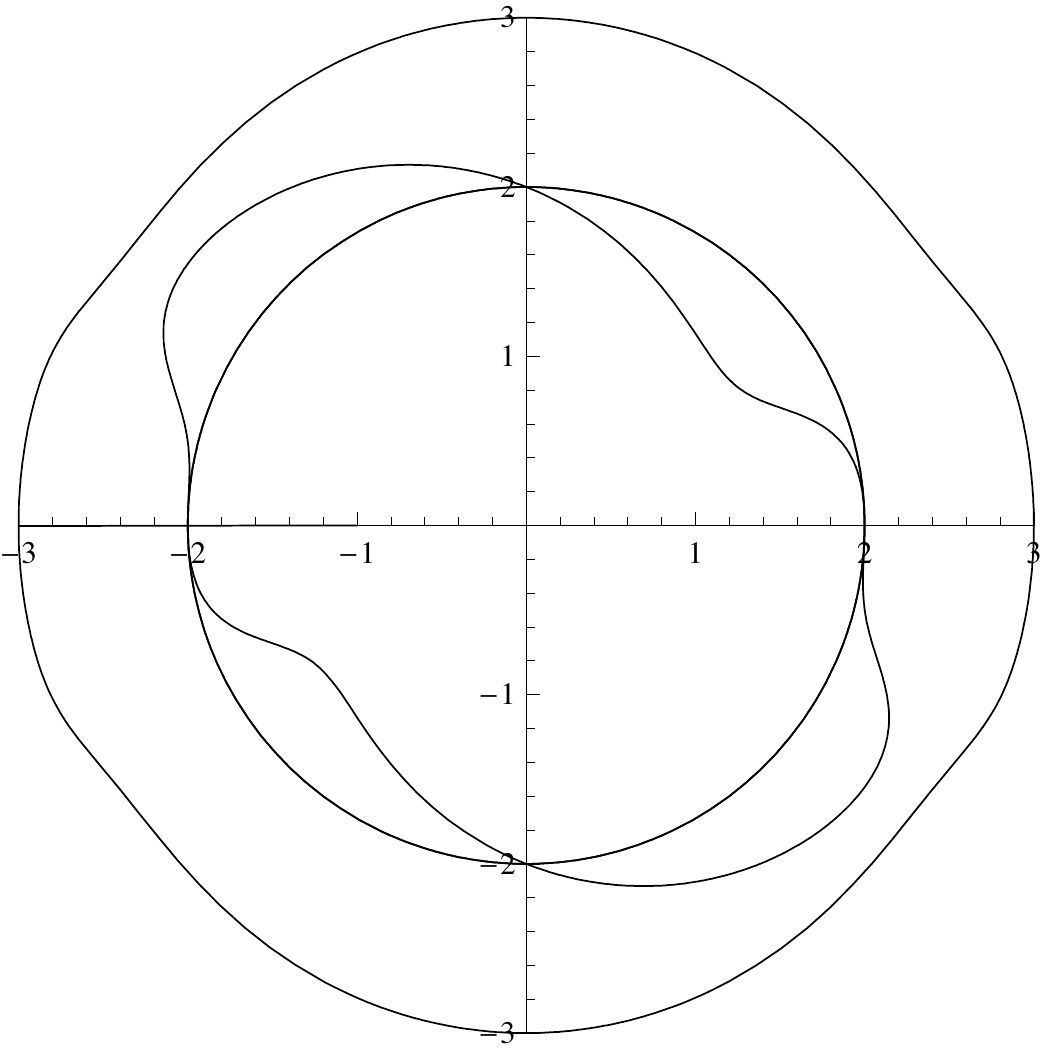}}
\caption{Fresnel curves and coupling functions for a rhombic medium with $\lambda=\frac12$, $\mu=1$,
$\tau_1=\frac52$ and $\tau_2=8$. The medium has four hyperbolic directions, in two the coupling function vanishes to third order.}
\label{fig2c}
\end{figure}

Except for the case with triple zeros, i.e. for $(\lambda+2\mu-\tau_1)(\lambda+2\mu-\tau_2)=0$, we are done and the results from \cite{RW07} imply the dispersive decay rate $t^{-1/2}$. In the exceptional
case we have to determine further the tangency index $\bar\gamma(\eta)$ of the corresponding
sheet of the Fresnel curve. Since for rhombic media the Fresnel curves are mirror symmetric with 
respect to the co-ordinate axes and the exceptional directions coincide with one of these axes, we
can only have $\bar\gamma_j(\eta)=2$ or $\bar\gamma_j(\eta)=4$ for the corresponding sheet. 
Because the hyperbolic eigenvalue in the exceptional direction is $\mu$ and the parabolic one is
$\tau=\lambda+2\mu$ this means that if $\lambda+\mu<0$ we are on the inner sheet and
$\bar\gamma_2(\eta)=2$. If $\lambda+\mu>0$ we are on the outer one and the higher order tangency
occurs according to \cite[Prop. 3.4]{RW07} if 
$\partial_\phi^2 \sqrt{\varkappa_1(\phi)}+\sqrt{\varkappa_1(\phi)}=0$
for the corresponding $\phi$. We consider the case $\tau_1=\lambda+2\mu$ such that the coupling
function $a_1(\eta)$ vanishes in $\phi=0$ to third order. To decide about the order of tangency it suffices to determine $\varkappa_1(\phi)$ up to $\mathcal O(\phi^3)$. 

We use $\eta=(\cos\phi,\sin\phi)^T=(1-\frac12 \phi^2,\phi)^T+\mathcal O(\phi^3)$, such
that
\begin{align}
&A(\phi) =\begin{pmatrix} (\lambda+\mu) (1-\phi^2)+\mu & (\lambda+\mu) \phi \\
 (\lambda+\mu) \phi & (\tau_2-\mu) \phi^2+\mu \end{pmatrix}+\mathcal O(\phi^3),\\
&r_1(\phi) = \begin{pmatrix} -\phi\\ 1-\frac12\phi^2\end{pmatrix} +\mathcal O(\phi^3),\\
&\varkappa_1(\phi) =  \mu + \beta \phi^2  +\mathcal O(\phi^3)
\end{align}
with an unknown constant $\beta$ to be determined via $A(\phi)r_1(\phi)=\varkappa_1(\phi)r_1(\phi)+\mathcal O(\phi^3)$. This implies 
\begin{equation}
A(\phi)r_1(\phi) = \begin{pmatrix} 
-\mu\phi   \\
\mu + (\tau_2-\lambda-\frac52\mu ) \phi^2\end{pmatrix} + \mathcal O(\phi^3)
\end{equation}
and
\begin{equation}
\varkappa_1(\phi)r_1(\phi)=\begin{pmatrix}
-\mu\phi  \\
\mu +(\beta-\frac12\mu) \phi^2 
\end{pmatrix}+\mathcal O(\phi^3)
\end{equation}
such that the unknown constants are determined by
\begin{equation}
    \beta-\frac12\mu = \tau_2-\lambda-\frac52\mu.
\end{equation}
This implies $\beta=\tau_2-\lambda-2\mu$. Thus $\varkappa_1(\phi) = 
\mu + (\tau_2-\lambda-2\mu) \phi^2 + \mathcal O(\phi^3)$ (even by symmetry of the Fresnel curve modulo $\mathcal O(\phi^4)$). We use that
\begin{align}
2\sqrt{\varkappa_1(\phi)} & \big( \partial_\phi^2 \sqrt{\varkappa_1(\phi)}
+\sqrt{\varkappa_1(\phi)} \big)\bigg|_{\phi=0}
=\partial_\phi^2\varkappa_1(\phi) - \frac{(\partial_\phi\varkappa_1(\phi))^2}{2\varkappa_1(\phi)}
+2 \varkappa_1(\phi) \bigg|_{\phi=0} \notag\\
&= \partial_\phi^2\varkappa_1(0) + 2\varkappa_1(0) = 2 (\tau_2-\lambda-2\mu) + 2\mu
= 2 ( \tau_2-\lambda-\mu). 
\end{align}
This expression vanishes if and only if $\tau_2=\lambda+\mu$. Hence, if $\tau_2\ne\lambda+\mu$
we get
$\bar\gamma_1(\eta)=2$ in this direction, while otherwise $\bar\gamma_1(\eta)=4$ holds true.

\begin{figure}[htbp]
\parbox{5cm}{\includegraphics[width=5cm]{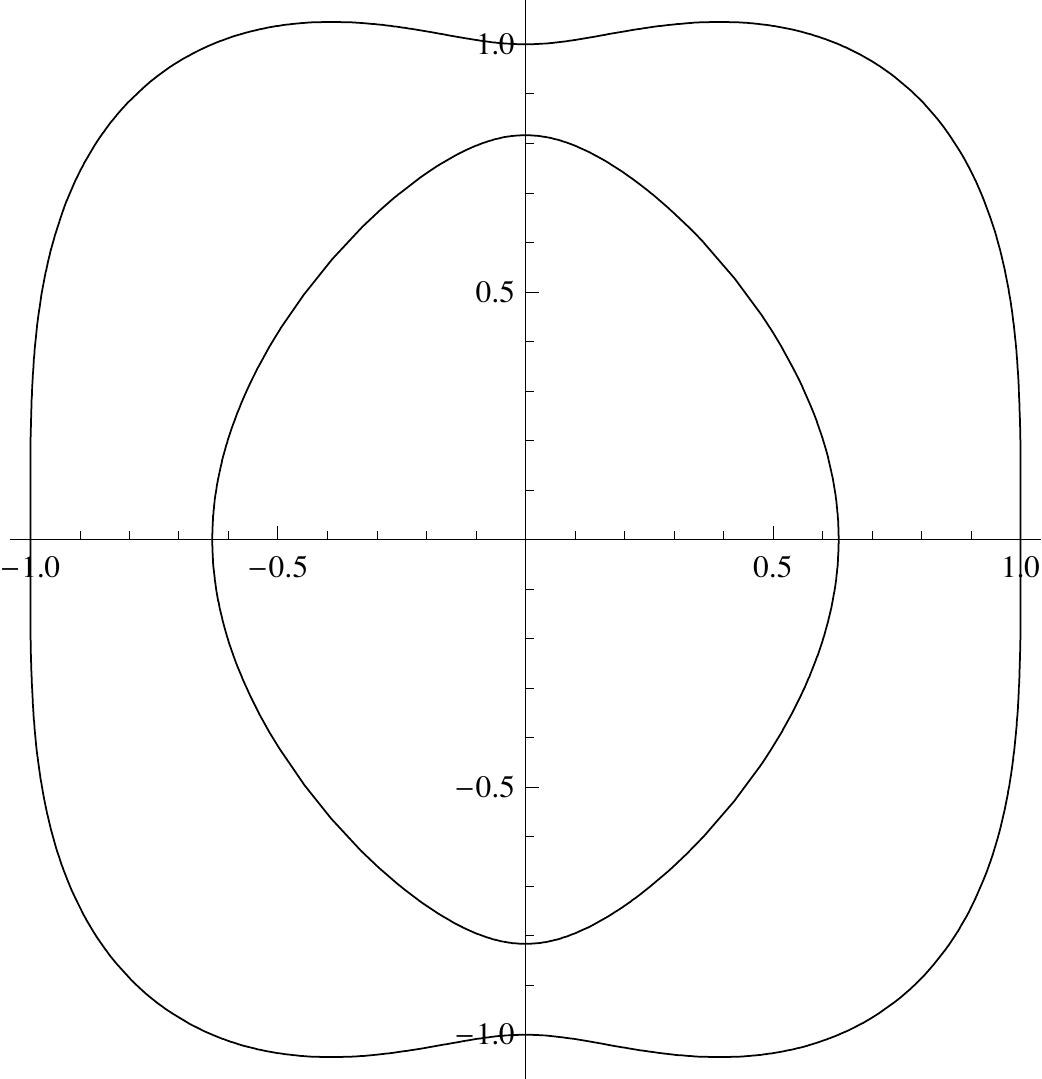}}
\parbox{5cm}{\includegraphics[width=5cm]{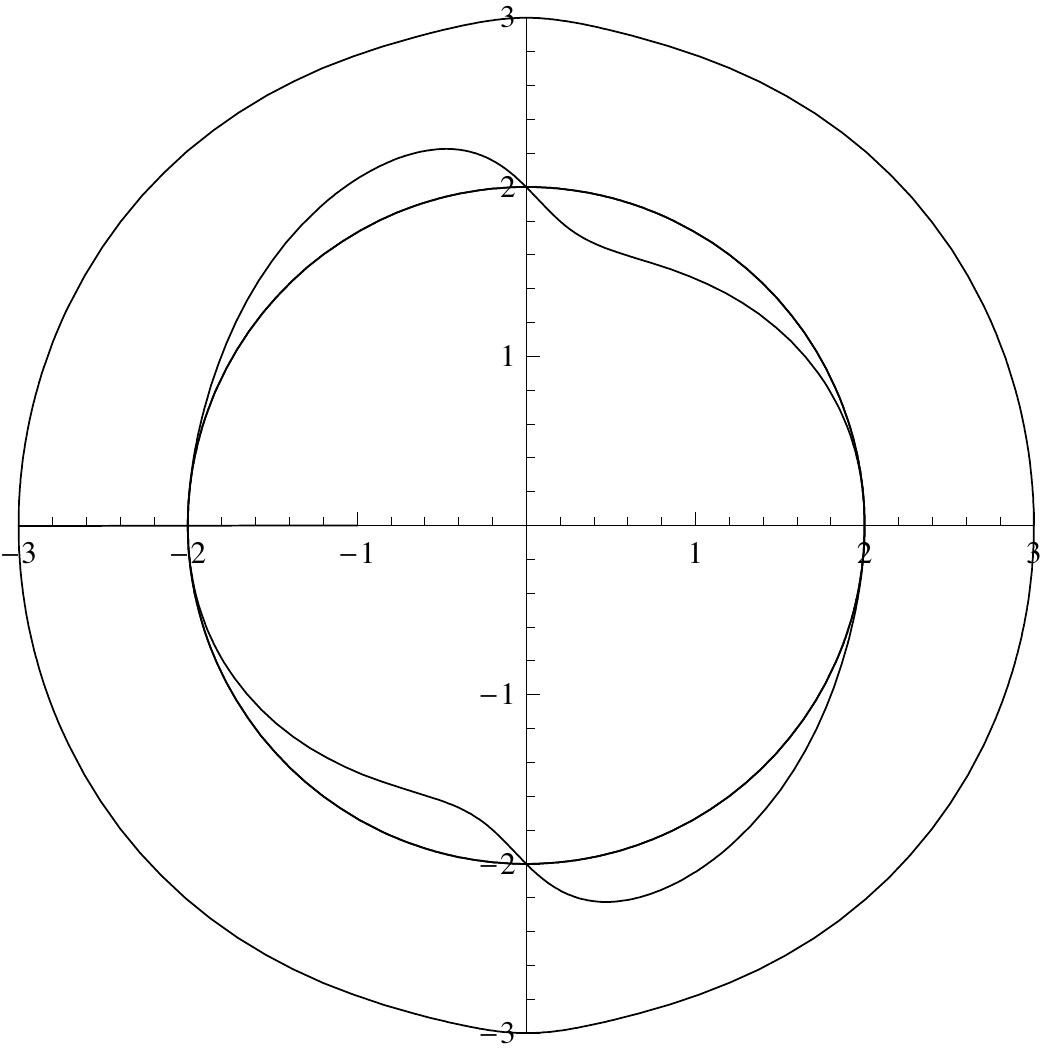}}
\caption{Fresnel curves and coupling functions for a rhombic medium with $\lambda=\frac12$, $\mu=1$,
$\tau_1=\frac52$ and $\tau_2=\frac32$. The medium has four hyperbolic directions, in two the coupling function vanishes to third order and the Fresnel curve has a higher order tangency.}
\label{fig2d}
\end{figure}

Finally we collect the conditions we have to impose for rhombic media:

(A1-2)\quad $\mu,\tau_1,\tau_2>0$, $-\sqrt{\tau_1\tau_2}-2\mu<\lambda< \sqrt{\tau_1\tau_2}$;

(A3)\quad $(\lambda+\mu)(\tau_1-\mu)(\tau_2-\mu)\ne0$;

(A4)\quad $\gamma^2\not\in\{\mu-\tau_1,\mu-\tau_2,-(\lambda+\mu)\}$ if $(\lambda+2\mu-\tau_1)(\lambda+2\mu-\tau_2)>0$ and \\
\hspace*{1.05cm}$\gamma^2\not\in\{\mu-\tau_1,\mu-\tau_2\}$ if $(\lambda+2\mu-\tau_1)(\lambda+2\mu-\tau_2)\le0$.

In Figures~\ref{fig2a} and \ref{fig2b} correspond to the situation of simple zeros of the coupling function, while in Figures~\ref{fig2c} and \ref{fig2d} one coupling function vanishes to order three.

\begin{cor} Rhombic media have the dispersive decay rate $t^{-1/2}$, except in the
case $\tau_1=\lambda+2\mu$ and $\tau_2=\lambda+\mu$ (or vice versa), where the dispersive decay
rate is $t^{-1/4}$. \end{cor} 

\section{An exceptional type of media}\label{sec6}
We want to consider media with $\tau_1=\tau_2=\lambda+2\mu$, $\sigma_1=0$ and $\sigma_2=\mu$. Again we start to check conditions (A1) to (A4). The matrix
$A(\eta)$ has the form
\begin{equation}
  A(\eta) = \begin{pmatrix} (\lambda+\mu)\eta_1^2+\mu & (\lambda+\mu)\eta_1\eta_2 + \mu \eta_2^2 \\ (\lambda+\mu)\eta_1\eta_2 + \mu\eta_2^2 &  (\lambda+\mu)\eta_2^2+\mu + 2\mu\eta_1\eta_2 \end{pmatrix}
\end{equation}
such that it is positive if and only if
\begin{equation}
   \mathrm{tr}\;A(\eta) = \lambda+3\mu + \mu\sin2\phi > 0,
\end{equation}
i.e. $\mu<\lambda+3\mu$, and
\begin{equation}
\det A(\eta) = \frac\mu8 (13\mu+8\lambda+4\mu\cos2\phi-\mu\cos4\phi+8\mu\sin2\phi+4(\lambda+\mu)\sin4\phi) >0.
\end{equation}
This is clearly positive if $\lambda>\mu>0$. [One can relax this condition to $\lambda>-\frac{8-3\sqrt3}4\mu$ and $\mu>0$. Note, that $-2+\frac34\sqrt3\approx-0.7$.]

Assumption (A3) is violated if and only if for some directions
\begin{equation}
  (\lambda+\mu)\eta_1\eta_2 + \mu\eta_2^2 =0\qquad\text{and}\qquad
   (\lambda+\mu)\eta_1^2=(\lambda+\mu)\eta_2^2+ 2\mu\eta_1\eta_2.
\end{equation}
The first equation implies $\eta_2=0$ or $(\lambda+\mu)\eta_1+\mu\eta_2=0$. In the first case
the second equation implies $\lambda+\mu=0$, while in the second one
the second equation implies $2(\lambda+\mu)\eta_1^2 = \lambda+\mu +2\mu\eta_1\eta_2$, which together with the first one gives $0=\lambda+\mu-4(\lambda+\mu)\eta_1^2$. Thus either 
$\lambda+\mu=0$ or $\eta_1^2=\frac14$. This corresponds to $\eta_2^2=\frac34$ and thus
we need $0=\lambda+\mu \pm \sqrt3\mu$, i.e. $\lambda + (1\pm \sqrt{3})\mu=0$.
The condition with $+$ can be skipped (since it is excluded by (A2)), the same for $\lambda+\mu$.  
Thus for (A3) we require $\lambda \ne (\sqrt{3}-1)\mu$.

Hyperbolic directions are determined by
\begin{equation}
0=  4\mu \cos2\phi-4\mu\cos4\phi,
\end{equation}
i.e. $\cos2\phi=\cos4\phi=2\cos^2(2\phi)-1$, which implies $\cos2\phi=\frac14(1\pm3)\in\{1,-\frac12\}$.
Thus the hyperbolic directions are $\phi=k\pi/3$ and we got double zeros for $\phi=k\pi$. We check condition (A4). For $\phi=0$ (or $\phi=\pi$) the eigenvalues of the matrix $A(\eta)$ are
$\mu$ (hyperbolic) and $\lambda+2\mu$, thus we exclude $\gamma^2\ne-(\lambda+\mu)$ (which
never happens since by (A2) $\lambda+\mu>0$). For $\phi=\pi/3$ (or $\phi=4\pi/3$) they are
$\mu-\sqrt3\mu/4$ (hyperbolic) and $\lambda+2\mu+3\sqrt3\mu/4$, thus we exclude
$\gamma^2\ne-\mu-\lambda-\sqrt3\mu$ (which also never happens due to (A2)). 
For $\phi=2\pi/3$ (or $\phi=5\pi/3$) we get
$\mu+\sqrt3\mu/4$ (hyperbolic) and $\lambda+2\mu-3\sqrt3\mu/4$, thus we have to exclude
$\gamma^2\ne-\lambda-\mu +\sqrt{3}\mu$.

It remains to consider the exceptional directions $\phi=k\pi$ and to investigate the vanishing order
of the coupling functions and the tangency index of the Fresnel curve. We have seen that the hyperbolic eigenvalue $\mu$ is always smaller than the other eigenvalue, thus we are on the outer sheet of
the Fresnel curve. Similar to the treatment of rhombic media we determine $\varkappa_1(\phi)$
modulo $\mathcal O(\phi^4)$. Because $\eta=(1-\frac12\phi^2,\phi-\frac16\phi^3)^T+\mathcal O(\phi^4)$
we get
\begin{align}
A(\phi) &= \begin{pmatrix} (\lambda+\mu)(1-\phi^2)+\mu & (\lambda+\mu) (\phi-\frac23\phi^3)+\mu \phi^2 \\ (\lambda+\mu) (\phi-\frac23\phi^3)+\mu \phi^2 & (\lambda+\mu)\phi^2 + \mu + 2\mu  (\phi-\frac23\phi^3)\end{pmatrix}+\mathcal O(\phi^4)\\
r_1(\phi) &= \begin{pmatrix}-\phi\\1 \end{pmatrix}+\phi^2\begin{pmatrix}\alpha_1\\
\alpha_2\end{pmatrix}+\phi^3\begin{pmatrix}\alpha_3\\
\alpha_4\end{pmatrix}+\mathcal O(\phi^4)\\
\varkappa_1(\phi)&=\mu+\beta_1\phi+\beta_2\phi^2+\beta_3\phi^3+\mathcal O(\phi^4)
\end{align}
with unknown constants  $\alpha_1, \alpha_2,\alpha_3,\alpha_4$ and $\beta_1, \beta_2, \beta_3$. Using
\begin{align}
  A(\phi)r_1(\phi) &= \begin{pmatrix}
 -\mu\phi+(\alpha_1(\lambda+2\mu)+\mu)\phi^2+(\alpha_3(\lambda+2\mu)+\alpha_2(\lambda+\mu)+\frac13(\lambda+\mu))\phi^3\\
 \mu+2\mu\phi+\alpha_2\mu\phi^2+(\alpha_4\mu+\alpha_1(\lambda+\mu) +2\alpha_2\mu-\frac73\mu)\phi^3
 \end{pmatrix}+\mathcal O(\phi^4),\\
  \varkappa_1(\phi)r_1(\phi) &=\begin{pmatrix}
   -\mu\phi+(\mu\alpha_1-\beta_1)\phi^2 +(\mu\alpha_3+\beta_1\alpha_1-\beta_2) \phi^3\\
   \mu+\beta_1\phi+(\mu\alpha_2+\beta_2)\phi^2+(\mu\alpha_4+\beta_1\alpha_2+\beta_3)\phi^3
  \end{pmatrix}+\mathcal O(\phi^4)
\end{align}
we see that $\beta_1=2\mu$, $\alpha_1=-3\mu/(\lambda+\mu)$, $\beta_2=0$, $\beta_3=-3\mu-7/3$ (while we would need an additional normalisation condition to determine $\alpha_2$, $\alpha_3$ and $\alpha_4$). Thus 
\begin{equation}
\varkappa_1(\phi) = \mu+2\mu\phi-(3\mu+\frac73)\phi^3+\mathcal O(\phi^4).
\end{equation}
This allows to determine $\bar\gamma_1$. Note first that
\begin{equation}
  2\sqrt{\varkappa_1(\phi)}(\partial_\phi^2\sqrt{\varkappa_1(\phi)}+\sqrt{\varkappa_1(\phi)}\big)\bigg|_{\phi=0}
  =\partial_\phi^2\varkappa_1(0)-\frac{(\partial_\phi\varkappa_1(0))^2}{2\varkappa_1(0)}+2\varkappa_1(0)=-2\mu+2\mu=0.
\end{equation}
Thus $\bar\gamma_1(\eta)>2$. On the contrary, 
\begin{align}
 2\sqrt{\varkappa_1(\phi)}\partial_\phi(\partial_\phi^2\sqrt{\varkappa_1(\phi)}+\sqrt{\varkappa_1(\phi)}\big)\bigg|_{\phi=0} &=
 \partial_\phi^3\varkappa_1(0)+(\partial_\phi \varkappa_1(0))\big(1-\frac{3\partial_\phi^2\varkappa_1(0)}{2\varkappa_1(0)}\big)+\frac{3(\partial_\phi\varkappa_1(0))^3}{4(\varkappa_1(0))^2}\notag\\
 &= -(18\mu+21)+2\mu+\frac{24\mu^3}{4\mu^2} = -10\mu-21 < 0
\end{align}
and therefore $\bar\gamma_1(\eta)=3$ for this direction.

Again we collect the conditions for this exceptional type of media

(A1-2)\quad $\lambda>\mu>0$;

(A3)\quad $\lambda\ne(\sqrt3-1)\mu$;

(A4)\quad $\gamma^2\ne(\sqrt3-1)\mu-\lambda$.  

\begin{figure}
\parbox{5cm}{\includegraphics[width=5cm]{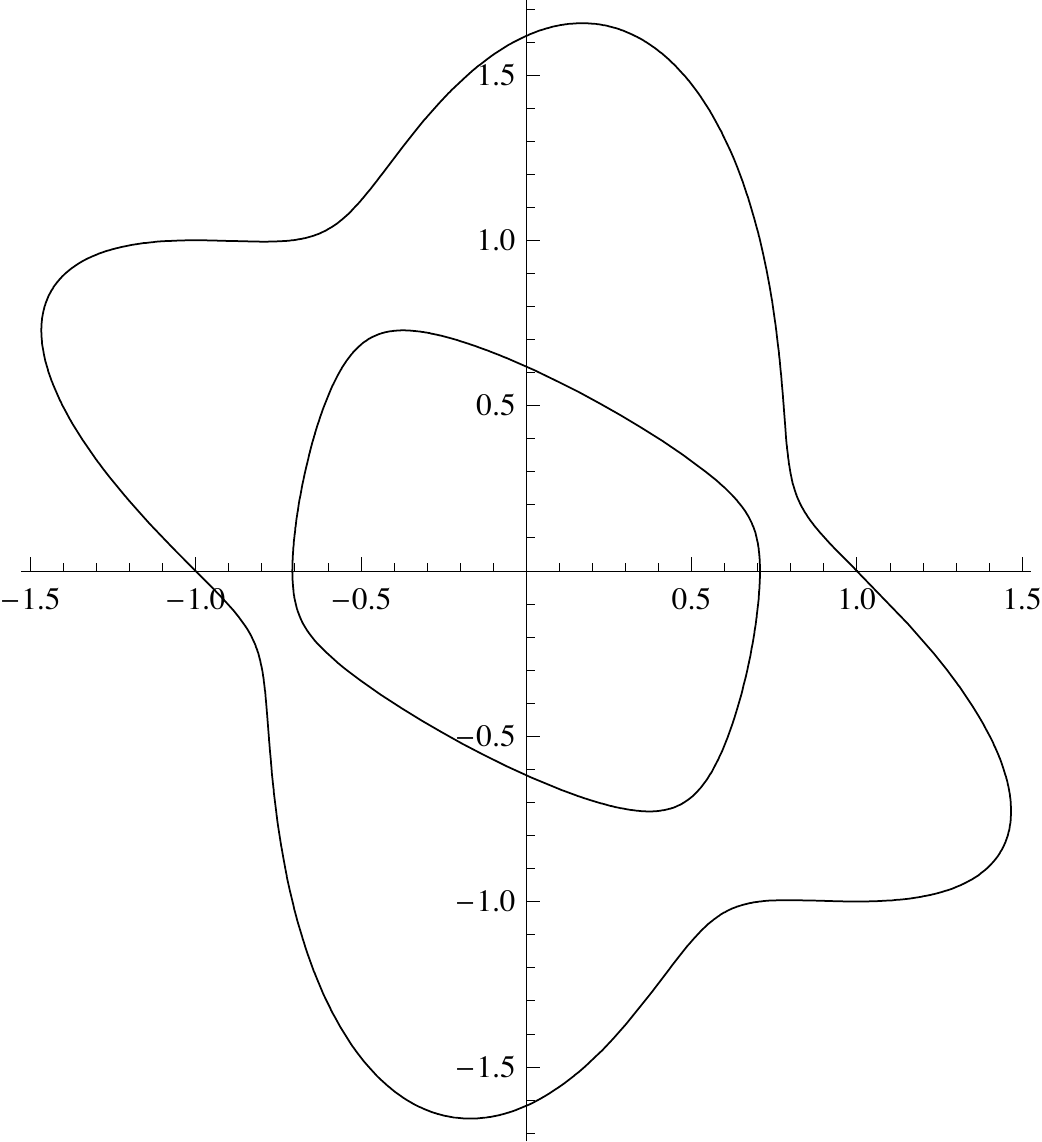}}
\parbox{5cm}{\includegraphics[width=5cm]{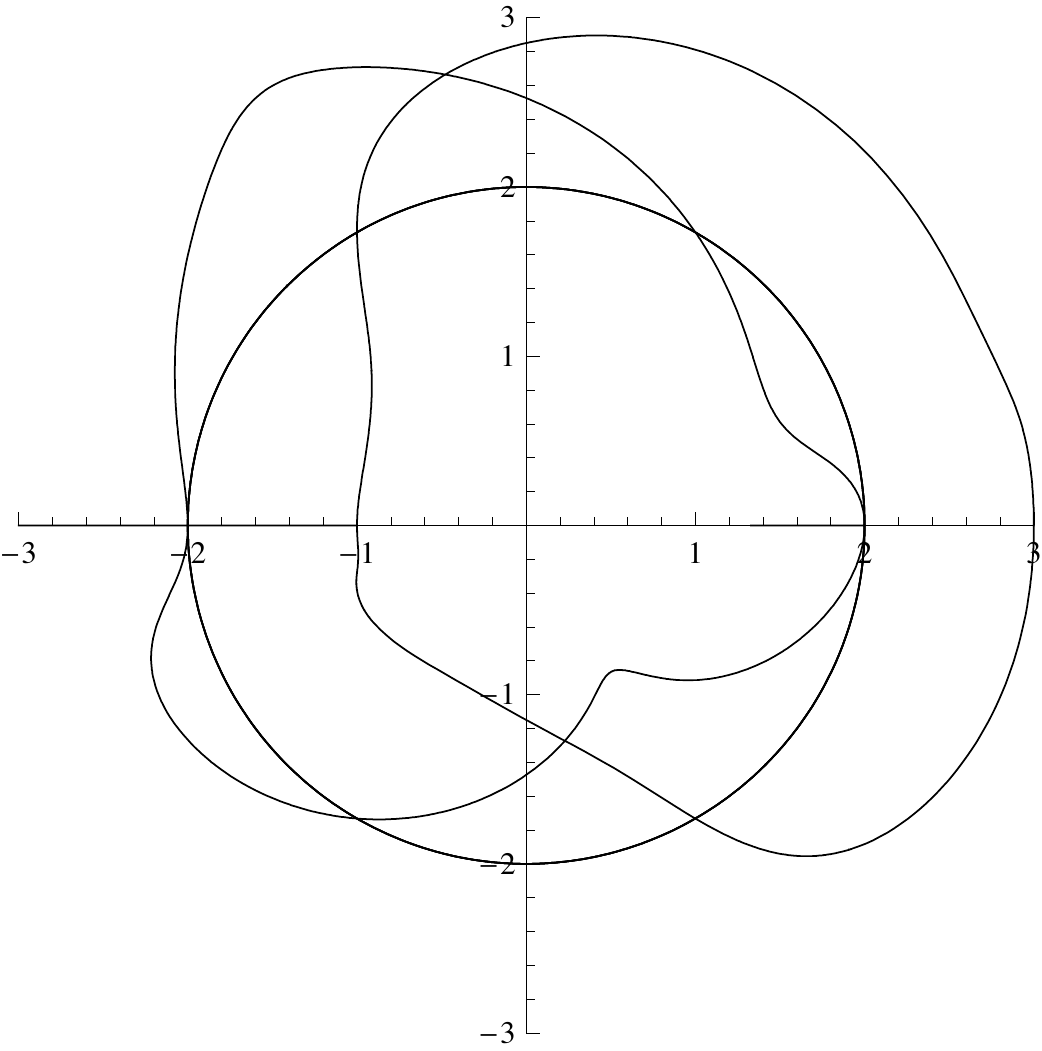}}
\caption{Fresnel curves and coupling functions for the exceptional medium with $\lambda=0$ and $\mu=1$. The medium has six hyperbolic directions, in two the coupling function vanishes to second order and the Fresnel curve has a point of inflection in these directions.}
\label{fig3}
\end{figure}

In Figure~\ref{fig3} the Fresnel curves and the coupling functions are depicted for one exceptional medium of this kind.

\begin{cor} This exceptional medium has the dispersive decay rate $t^{-1/3}$. \end{cor}

\section*{Acknowledgements}
The author thanks the Department of Mathematics, University College London, for hospitality and support during the last year.
The figures have been produced using the computer algebra system Mathematica 6.0. The author is grateful to the University College London for providing these facilities. 


\end{document}